\documentclass[12pt]{article}
\usepackage{amsmath}
\usepackage{amsbsy}
\usepackage{amssymb}
\usepackage{a4wide}
\newenvironment{proof}{\par\vspace{1ex}\par {\em Proof:}\hspace{0.5em}}{\noindent$\bigtriangleup$\par\vspace{1ex}\par }
\newtheorem{corollary}{Corollary}

\newtheorem{theorem}{Theorem}

\newtheorem{definition}{Definition}
\newtheorem{example}{Example}

\begin{document}

\title{Gr\"obner bases of $*$-algebras and faithful operator representations.}

\author{S. Popovych, K. Vynogradov}

\maketitle

\newcommand{\ba}[1]{\ensuremath{e_{\phi(#1)}}}
\newcommand{\sca}[2]{\ensuremath{\langle #1,#2 \rangle}}
\newcommand{\algebra}[1]{\ensuremath{\mathbb{C}\langle #1 \rangle}}

\begin{abstract}
We give  sufficient conditions for existence of a faithful
representation of a $*$-algebra in terms of its G\"obner basis. In
order to do this we propose a construction of a faithful
representation. This construction applies to concrete examples:
$*$-doubles, monomial $*$-algebras, extension of a $*$-algebras
allowing Wick ordering and others. Several examples and
counterexamples are also presented.
\end{abstract}


\section{Introduction}
\label{intro}

\noindent

We address a fundamental question in representation theory: When
can an algebra (with involution) built on generators and relations
be represented by unbounded operators on Hilbert space? And when
is there a faithful representation? We make a case for the
usefulness of Gr\"obner basis techniques.   Now, it has been known
since the early days of the algebraic approach to quantum problems
(e.g., Heisenberg's commutation relations) that there are
representation theoretic dichotomies: For example, the relations
of the Bosons and the Fermions display different features when
represented in Hilbert space. This is the issue of unbounded
operators vs bounded operators. Both examples fit the theme of the
paper, viz., algebras A with involution over the complex field. In
the special cases when $A$ is assumed abelian, or if $A$ is a
$C^*$ -algebra, we have the familiar theorems of Gelfand and
Naimark, but for algebras on generators and relations, the
literature so far only consists of isolated classes of examples.
In particular, the question of representability by Hilbert space
operators  has been studied for $*$-algebras allowing Wick
ordering, monomial $*$-algebras and other
(see~\cite{tapper,lance,jorgensen,prosk1,pop1}).  Lance and
Tapper~\cite{tapper,lance} undertook an attempt to treat
$C^*$-representability of monomial $*$-algebras with one defining
relation. However there isn't a general theory yet.

Let us fix some notations.  We will denote by  $E$ a pre-Hilbert
space and by $H$ a Hilbert space. Let  $L(E)$ and $L(H)$ denote
the $*$-algebras of linear operators acting on these spaces. We
will study  the question whether  a $*$-algebra $A$ can be
embedded as a $*$-subalgebra in $L(E)$ or $L(H)$. Note that in the
latter case elements of $A$ are represented  by bounded operators.
 If a $*$-algebra  $A$ is embedded in $L(E)$ and every operator
$a\in A$ is bounded then one can extend each  $a\in A$ to an
operator acting on the completion  $H$ of $E$ and thus obtain an
inclusion $A\hookrightarrow L(H)$. In the  general case $A$ will
be represented by unbounded   operators  on $H$ such that the
intersection of their  domains is dense.  Henceforth when
unbounded operators are involved  the term {\it representation }
of an algebra means that the relations of the algebra are
satisfied on the common dense invariant  domain.
 Remind that $*$-algebra  $A$ is called $C^*$-representable if
 there is its faithful $*$-representation by bounded operators  acting on
 a Hilbert space or, equivalently, $A$ can be embedded into
 $C^*$-algebra.   Thus we will decompose  the question of $C^*$-representability
 into two parts: the first, of algebraic nature, is to find a
faithful representation of $A$ in a pre-Hilbert space and the
second  one, of analytic nature, is to find out whether this
representation is in bounded operators.

The question of $C^*$-representability  of finitely-presented  $*$-algebras
 has been investigated in~\cite{pop1}. The main novelty of
 the approach used in~\cite{pop1} was
  the employment of Gr\"obner basis
 technique.  Here we will elaborate this approach.

In present work we find sufficient conditions for a $*$-algebra to
be faithfull represented in pre-Hilbert spaces.  These conditions
expressed in terms of its defining relations and can be
effectively verified.

For $*$-algebra given by generators and relations if not the only
then at least very natural way to prove  that a homomorphism is
injective is to show that some linear  basis is mapped into a
linear basis of the image. To  construct a linear basis we use the
machinery of Gr\"obner bases developed in~\cite{Bok}. A brief
account of it  is given in the appendix.

The main property of a $*$-algebra which enable us to construct a
faithful representation is the strongly non-expanding condition
(see definition~\ref{nonexpand}). But this property is hard to
verify in examples. So we gave several sufficient conditions (see
definition~\ref{approp} and theorem~\ref{kir}) and apply them to
many concrete examples: $*$-doubles,  monomial $*$-algebras,
extension of a $*$-algebras with Wick ordering and others.

\section{ Unsrinkability type restrictions on Gr\"obner basis of a *-algebra.}
We will denoted by $F$ the free associative algebra since the
number of generators is not important for further considerations.
 In this paper we will deal exclusively with  finitely-generated algebras. Let us denote by  $F_*$ the
  free associative algebra with generators $x_1,x_2,\ldots,
x_m,x^*_1,x^*_2,\ldots, x^*_m$. Algebra  $F_*$ is a $*$-algebra
with involution given on generators by $(x_j)^*=x_j^*$ for all $j
=1,\ldots, m$. This is just associative algebra with $2m$
generators $F_{2m}$ with a natural involution. Note that $F_*$ is
a semigroup algebra of a semigroup $W$ of all words in generators
$x_1,x_2,\ldots, x_m,x^*_1,x^*_2,\ldots, x^*_m$. A set $S\subseteq
F$ of defining relations of an associative algebra is called
Gr\"obner basis if it is closed under compositions (see Appendix).
 A Gr\"obner basis
of a $*$-algebra $A$ is a Gr\"obner basis of $A$ considered as an
associative algebra. We need to put some extra requirements on a
Gr\"obner basis to make it "compatible" with involution. The main
requirement we impose is a generalization of the notion of
unshrinkability of the word.  A set $S\subseteq F_*$ is called
{\it symmetric} if the ideal $\mathcal{I}$ generated by $S$ in
$F_*$ is a $*$-subalgebra of $F_*$. In particular, $S$ is
symmetric if $S^*=S$.
  The following definition is due to P.Tapper~\cite{tapper}:
\begin{definition}
The word $w\in{W}$ is called unshrinkable if it
 can not be presented in the  form $w=d^*du$ or $w=ud^*d$
 for some nonempty word d.
\end{definition}

Recall that a $*$-algebra $A$ is  called {\it proper} if for every
element $x \in A$ condition $x^*x=0$ implies $x=0$. A $*$-algebra
$A$ is called {\it completely proper} if $M_n(\mathbb{C})\otimes
A$ is proper for all integer $n$. The importance of this notion
follows from the fact that any bounded unital simple $*$-algebra
is $C^*$-representable if and only if it is completely
proper~\cite{pop1}.

 P.Tapper has conjectured that $*$-algebra $\algebra{x,x^*|w=0, w^*=0}$
 is $C^*$-representable if and only if  word $w$ is unshrinkable.  The
first author in \cite{pop2} proved that monomial $*$-algebra $A$
is completely proper if and only if its defining relations are
unshrinkable. Moreover, in this case $A$ can be faithfully
represented by operators acting on pre-Hilbert space. A much wider
class of $*$-algebras for which we will prove similar results is
defined below. For the notations $u\prec w$, $R_S(w)$, $BW$ and
order on $W$ used below we refer the reader to the appendix.
\begin{definition}\label{nonexpand}
A symmetric subset  $S\subseteq F_*$ closed under compositions
will be called { non-expanding} if for every $u, v, w \in BW $
such that $u\not= v$ and $ww^*\prec {\rm R}_{S}(uv^*)$ the
following inequality holds $w< \sup{(u,v)}$. If in addition for
every word $d \in BW$ the word $dd^*$ also belong to  $BW$ we will
call $S$ { strictly non-expanding}. A $*$-algebra $A$ is called
non-expanding
 if it possesses a Gr\"obner basis $GB$ which is
 non-expanding and $A$ is strictly non-expanding if some of its Gr\"obner bases is strictly non-expanding.
\end{definition}
We will show (see theorem~\ref{main}) that every strictly
non-expanding algebra have a faithful representation in a
pre-Hilbert space.

\section{Representation construction.}

In this section we will show that strictly non-expanding
$*$-algebra possesses  a faithful positive functional and thus  is
a pre-Hilbert $*$-algebra.

 Let $G \subseteq W_n$  be a subset. An
 {\it enumeration} of $G$  is  a map $\phi:G \rightarrow
\mathbb{N}$ such that $u > v$ implies $\phi(u)>\phi(v)$ which is
bijection on $\mathbb{N}$ if $|G|=\infty$ and on some interval
$[1,n]\subseteq \mathbb{N}$ if $|G|<\infty$.

\begin{definition}\label{form}
Let us define the Half-word operator ${\rm H}:F_*\to F_*$ by the
rule ${\rm H}(uu^*)=u$ for $u\in W$ and ${\rm H}(v)=0$ if $v$ is
not positive word and extend by linearity to $F_*$. Let us fix
 a set $S\subseteq F$ closed under compositions, an enumeration
$\phi:BW\to\mathbb{N}$ of the corresponding linear basis and a
sequence of positive real numbers $\xi=\{a_k\}_{k \in \mathbb{N}}
\subset \mathbb{R}$. Define a weight functional
$T_{\xi}^{\phi}:K\to \mathbb{C}$ by putting
$T_{\xi}^{\phi}(u)=a_{\phi(u)}$ for any word $u\in BW$, where $K$
denote the linear span of $BW$. Let $n=|BW|\in \mathbb{N}\cup
\{\infty\}$.

 Let $V$ be a vector space
 over $\mathbb{C}$ with basis
 $\{e_k\}_{k = 1}^n$. Let us define a form $\sca{e_i}{e_j}_\xi$ on
 the basis elements and extend it by linearity in the first
 argument and conjugate-linearity in the second one by the following rules. If
$i=j$ then $\sca{e_i}{e_i}_\xi = a_i$. If $i\not=j$ then there are
unique elements $u,v\in BW$
  such that $i=\phi(u)$, $j=\phi(v)$. Put
  $\langle e_i,e_j \rangle_\xi = T_{\xi}^{\phi}\circ
{\rm H} \circ {\rm R}_{S}(uv^*)$.

\end{definition}

\begin{theorem}\label{scal}
 If $S$ is strictly non-expanding then there exist sequence
$\xi_0=\{ a_k \}_{k \in \mathbb{N}} \subset \mathbb{N}$ such that
the sesquilinear form $\langle \cdot , \cdot \rangle_{\xi_0}$ is
positively defined.
\end{theorem}
\begin{proof}
That the form $\langle \cdot , \cdot \rangle_\xi$ is sesquilinear
is obvious by definition.  Let $g_{ij}=\langle e_i,e_j
\rangle_\xi$ for $i,\ j\in \mathbb{N}$ and $G=(g_{ij})_{1\le
i,j\le\infty}$ denote the  Gram matrix of the sesquilinear form.
We will use Silvester's criterion and show
 by induction on $m$ that $a_m$ can be chosen such that principal minor $\Delta_m
 > 0$. For $m=1$ put $a_1=1$ then $\Delta_1 = 1 > 0$. Let $a_1,\ldots,a_{m-1}$ be already
chosen such that $\Delta_1>0,\ldots, \Delta_{m-1} > 0$.

 By definition if $u \in BW$, then $u^*u$ is  also
in $BW$. Thus by  definition~\ref{form}  we have $\langle
e_{\phi(u)},e_{\phi(u)} \rangle_\xi = a_{\phi(u)}$. Take some
$i\le m$ and $j\le m$  with $i\not=j$ and find unique $u, v\in BW$
such that $i=\phi(u),\ j=\phi(v)$. Then $uv^*=\sum_{k}\alpha_k
w_k$ for unique $\alpha_k\in \mathbb{C}$ and $w_k\in BW$. Clearly
$\langle e_{\phi(u)},e_{\phi(u)} \rangle_\xi$ is a sum $\sum_k
\alpha_k a_{\phi(h_k)}$ where the sum is taken over those $k$ for
which $w_k$ is a positive word $w_k=h_kh_k^*$. By  definition we
get that $h_k<\sup{(u,v)}$. Hence $g_{ij}$ is a polynomial in
variables $a_1,\ldots,a_{m-1}$. Now from the decomposition by the
$m$-th row we obtain $\Delta_m=\Delta_{m-1}a_m +p_{m}
(a_1,\ldots,a_{m-1})$, where $p_{m} \in \mathbb{C}[a_1,\ldots,
a_{m-1}]$ some polynomial. Since $\Delta_{m-1}>0$ it is clear that
$a_m$ can be chosen such that $\Delta_m>0$. This
 completes inductive proof.
\end{proof}

The space $K$ is isomorphic to $V$ via the map $u\to e_{\phi(u)}$.
Thus the inner product $\sca{\cdot}{\cdot}_\xi$ on $V$ give rise
to an inner product on $K$ which will be denoted by the same
symbol. It is routine to check that
$\sca{u}{v}_\xi=\alpha(P(u\diamond v^*))$ where $P:F\to F$  is the
projection on the linear span of positive words $W_+$ and $\alpha:
K\to \mathbb{C}$ some linear functional. Let $z\to L_z$ denote the
right regular representation of $A=F/\mathcal{I}$,  $L_z(f)=f z$
for any $z,f \in A$.
\begin{theorem}\label{main}
Let $S\subseteq F$ be strictly non-expanding and $\mathcal{I}$ the
ideal generated by $S$ in $F_*$. Then the right regular
representation $L$ of the $*$-algebra $A=F_*/\mathcal{I}$ on a
pre-Hilbert space $(K,\sca{\cdot}{\cdot}_\xi)$ is a faithful
$*$-representation.
\end{theorem}
\begin{proof}
The representation stated in the theorem is associated by the GNS
construction with the positive functional $\alpha(P(\cdot))$ on
$A$. Thus it is a $*$-representation. Indeed, as in the  GNS
construction $N=\{a \in A | \alpha(P(aa^*))=0\}$ is a right ideal
in $A$. We can define inner product on $A/N$ by usual rule
$\sca{a+N}{b+N}=\alpha(P(a^*b))$. It is easy to verify that right
multiplication define $*$-representation of $A$ on pre-Hilbert
$A/N$. The only difference with classical  GNS construction is
that this representation could not be, in general, extended on
completion of $A/N$.

Let us show that $*$-representation is a faithful
$*$-representation. Take any $f={\sum_{1 \le i \le n} c_i w_i}\in
A$, where $c_i \in \mathbb{C}, w_i \in BW$. Without loss of
generality consider $w_1$ to be the greatest word among $w_j$.
Then $L_f(w_1^*)$ contains element $w_1w_1^*$ with coefficient
$c_1$. Hence $L_f\not=0$.
\end{proof}
\begin{corollary}
Every strictly non-expanding $*$-algebra has a faithful
*-representation by unbounded operators.
\end{corollary}

\section{Sufficient conditions. Examples.}

We  call a subset $S\subseteq F$ {\it  reduced} if for any $s\in
S$ and any word $w\prec s$ no word $\hat{s'}$ with $s'\in S$ is
contained in $w$ as a subword. If the set $S$ is closed under
compositions then one can obtain reduced set $S'$ closed under
compositions
 generating the same ideal by replacing each $s\in S$ with
 ${\rm R}_S(s)$.
The following, rather technical,  modification of the notion of
{\it appropriate} $*$-algebra from~\cite{pop1} is the main source
of examples.
\begin{definition}\label{approp}
A symmetric reduced  subset $S\subseteq F_*$ is called  strictly
appropriate if it is closed under compositions and for every $s
\in S$ and every word $u\prec s$ such that $|u|=deg(s)$ the
following conditions hold:
\begin{enumerate}
\item   word $u$ is unshrinkable; \item  if  $u\not=\hat{s}$,
$\hat{s}=a b$, and $u=a c$ for some words $a,b,c$ with
$a\not=\emptyset$ then for any $s_1\in S$ such that there is word
$w\prec s_1,\ w\not= \hat{s}_1, |w|=|\hat{s}_1|$  either word
$\hat{s}_1$ does not contain $u$ as a  subword  or $\hat{s}_1$ and
$u$ do not form a composition in such a way that $\hat{s}_1=d_1 a
d_2$ and $u=a d_2 d_3$ with some nonempty words $d_1, d_2, d_3$.
\end{enumerate}
 A $*$-algebra $A$ is called  strictly appropriate  $*$-algebra if it possesses a strictly appropriate Gr\"obner
basis.
\end{definition}

The following lemma shows that strictly appropriate $*$-algebras
provide examples of non-expanding $*$-algebras. Many concrete
examples of strictly appropriate $*$-algebras will be considered
in the final section. In the following lemma for word $w\in W$ of
even length $w=w_1w_2, |w_1|=|w_2|$ we will denote by ${\rm
H}_0(w)$ the first half of $w$,   ${\rm H}_0(w)=w_1$.
\begin{theorem}\label{principal}
Every strictly appropriate set $S\subseteq F$ is non-expanding. If
in addition  $S=S^*$ then $S$ is strictly non-expanding.
\end{theorem}
\begin{proof}
Let $u, v\in BW$ and $u\not=v$.

1. If $uv^*\in BW$ then $ww^*\prec {\rm R}_{S}(uv^*)$ implies
$ww^*=uv^*$. By lemma 2~\cite{pop1} two cases are possible (1)
$u=vdd^*$ and $w= v d$ or (2) $v=udd^*$ and $w = v d$ where word
$d$ is nonempty. Hence $|w|=|d| +|v|= |u|-|d^*|< |u|$ in the first
case and  $|w|=|d| +|u|= |v|-|d^*|< |v|$ in the second one. Thus
$w< u$ or $w<v$.

2. Now let $uv^*\not\in BW$. There are words $p,q \in BW$ and
element $s\in S$ such that $uv^*=p\hat{s}q$. Moreover, since $u,v
\in BW$ none of them can contain  $\hat{s}$ as a subword. Hence
$\hat{s}=a b$ with nonempty words $a$ and $b$ such that $u=p a$
and $v^*= b q$. Write down $s=\alpha \hat{s} + \sum_{i=1}^k w_i
+f$, where $deg(f) < deg(s)$ and $|\hat{s}|=|w_i|$ for all $i\in
\{1,\ldots, k\}$. Assume that for some integer $i$ word $pw_iq$
belongs to $BW$ and $pw_iq=ww^*$ for some word $w$. If the middle
of the word $pw_iq$ comes across $w_i$, in other words
$max(|p|,|q|)<|w|$ then $w_i=c d$,  $w= p c$,  and $w^*=d q$ with
some nonempty words $c,d$. Hence $p c= q^*d^*$. If $|c|\le |d|$
then $d^*=g c$ for some word $g$ and so $w_i= c d = c c^* g^*$
which contradicts unshrinkability of $w_i$. If $|c|>|d|$ then $p
c= q^*d^*$ implies $c= g d^*$ for some word $g$ and we again see
that $w_i= g d^*d$ is shrinkable. Thus $max(|p|,|q|)\ge |w|$. If
$|p|>|w|$ then $|u|=|p|+|a|>|w|$, otherwise $|v|=|b|+|q|>|w|$.

In the above-mentioned  cases we have had $|w|<max(|u|,|v|)$ which
is a stronger statement than that of the lemma. But on the second
step of the decomposition process which we now approaching this
regularity breaks down.

3. Let $uv^*=p\hat{s}q$ and $s=\alpha \hat{s} + \sum_{i} w_i +f$
as above and $ww^*\prec {\rm R}_{S}(pw_iq)$. It is obvious that
$w\le \sup{(u,v)}$. We need to prove that $w< \sup{(u,v)}$.
Suppose the contrary: $uu^*\prec {\rm R}_{S}(uv^*)$ and $u>v$ or
$vv^*\prec {\rm R}_{S}(uv^*)$ and $v>u$. If $|u|\not=|v|$ then
$ww^*\prec R_{S}(pw_iq)$  implies $|w|<max(|u|,|v|)$ because $d<c$
for all  $d\prec {\rm R}_{S}(c)$. If $|\hat{s}|>|w_i|$ then again
$|w|<max(|u|,|v|)$. Hence we can assume that $|u|=|v|$ and
$|\hat{s}|=|w_i|$. Clearly $u<v$ implies $uv^*<vv^*$. Hence
relation $v v^* \prec {\rm R}_S(p w_i q)$ is impossible. Thus we
are left with the only possibility $u>v$ and $uu^*\prec {\rm
R}_{S}(p w_i q)$. Since $uu^*<p w_i q<uv^*$ word  $p w_i q$ begins
with $u$. If $\hat{s}=a b$ such that $p a=u, b q =v^*$ then $w_i$
should also begin with $a$. Therefore  $\hat{s}$ and $w_i$ begin
with the same generator. Since $pw_i q\not\in BW$ there is
$s_1=\alpha \hat{s}_1 +\sum_j\beta_j u_j +g\in S$ with
$deg(g)<deg(s_1)$ such that $pw_iq=p_1\hat{s}_1q_1$ for some words
$p_1, q_1$. If we assume that  for some $j$ word $uu^*\prec {\rm
R}_{S}(p_1u_jq_1)$ then $|\hat{s}_1|=|u_j|$ and ${\rm
H}_0(p_1u_jq_1)=u$. Word $\hat{s}_1$ could not be a subword in the
first half of the word $pw_i q$ since ${\rm H}_0(p_1u_j q_1)={\rm
H}_0(pw_iq)=u$ and assuming the contrary we see that $\hat{s}_1$
and $u_j$ are both subwords of $u$ in the same position, hence
should be equal $\hat{s}_1=u_j$.
 Word $\hat{s}_1$ could not contain subword $w_i$ because of condition $1$
 in the definition of strictly appropriateness.
 Obviously, $\hat{s}_1$ could not be a subword in $q$ because  $q\in
BW$. Thus either $w_j$ and $\hat{s}_1$  intersect or $\hat{s}_1$
and  $w_j$ intersect in such a way that $\hat{s}_1=d_1 a d_2$ and
$w_j = a d_2 d_3$. But this contradicts the strictly
appropriateness of $S$. So we have proved that $w<\sup{(u,v)}$.
The fact that for any word $g\in BW$ word $gg^*$ lies in $BW$
follows from lemma 2.~\cite{pop1}.
\end{proof}

The following is a convenient simplification of the preceding
theorem:
\begin{corollary}\label{cor}
If a symmetric subset $S\subseteq F$ is  closed under
compositions,
 for every $s \in S$ and every word $u\prec s$ such that $|u|=deg(s)$ the word $u$ is unshrinkable and
  words $\hat{s}$ and $u$  begin with different generators then $S$ is non-expanding.
  If in addition $S=S^*$ then $S$ is strictly non-expanding.
\end{corollary}

\begin{example} Let $\mathcal{L}$ be a finite dimensional real Lie
algebra with linear basis $\{e_j\}_{j=1}^{n}$. Then its universal
enveloping algebra $U(\mathcal{L})$ (over $\mathbb{C}$) is a
$*$-algebra with involution given on generators as $e_j^*=-e_j$.
We claim that this $*$-algebra is non-expanding.
  Indeed $M=\{ e_i
e_j-e_j e_i- [e_i,e_j], i<j \}$ is a set of defining relations for
$U(\mathcal{L})$ it is closed under compositions (see example
in~\cite{Bok} or use PBW theorem).  Then the set $S= \{e^*_j+e_j,
1\le j \le n \}\cup M$ is also closed under compositions (we
consider $e^*_1>e^*_2>\ldots>e^*_1>e_1>\ldots > e_n$) since
$e^*_j$ and $e_k e_l$ do not intersect for any $j$, $k$, $l$. It
is easy to see that $S$ is symmetric. Thus $S$ is non-expanding by
corollary\ref{cor}.
\end{example}

\begin{theorem}\label{kir}
Let $S\subseteq F_*$ be a symmetric subset of a free countably
generated $*$-algebra $F_*=F(X\cup X^*)$ such that the following
conditions are satisfied:
\begin{enumerate}
\item For every  $s\in S$ and every word $w\prec s$ with
$|w|=deg(s)$ is unshrinkable. \item For every $s_1$, $s_2\in S$
and every word $u\prec s_1$ with $|u|=deg(s_1)$ the words $u$ and
$\hat{s_2}$ do not form a composition. Then  $*$-algebra
$A=\mathbb{C}\langle X\cup X^*| S\rangle$ is non-expanding. If in
addition $S=S^*$ then $A$ is strictly non-expanding.
\end{enumerate}
\end{theorem}
\begin{proof}
It is suffices to show that for any two basis words $u$, $v\in BW$
such that $|u|=|v|$ and $u>v$ and for any word $p$ condition
$pp^*\prec {\rm R}_S(uv^*)$ implies $pp^*<uu^*$.

Indeed, if $|u|\not= |v|$ then, clearly, $|pp^*|\le |uv^*| <
max(|uu^*|,|vv^*|)$. Hence $p <\sup{(u,v)}$. If $|u|=|v|$ and
$v>u$ then $uv^*<vv^*$. Thus  every word $w$ such that $w\prec
{\rm R}_S(uv^*)$ is less than $vv^*$.

 Since $pp^*\le uv^*$ it remains to prove
that $uu^*\not\prec {\rm R}_S(uv^*)$. Assume the contrary. Then
there is a sequence of words $\{q_i\}_{i=1}^{n}$ such that
$q_1=uv^*$, $q_n=uu^*$ and for every $1\le i\le n-1$ there is
$s_i\in S$ and words $c_i$, $d_i$, $u_i\in W$ such that $u_i\prec
s_i$, $u_i\not=\hat{s_i}$, $|u_i|=|\hat{s}_i|$ and $q_i=c_i
\hat{s}_i d_i$, $q_{i+1}=c_i u_i d_i$.

Let $j$ be the greatest with  the property that $\hat{s}_j$
intersects the middle of $q_j$. Such an index $j$ exists because
$j=1$ satisfies this property and we choose within a finite set.
Clearly $j<n$ since otherwise $u_{n-1}$ would be a subword in
$uu^*$ intersecting its middle and thus would be shrinkable. Thus
for every $i\in \{j+1,\ldots,n-1 \}$ word $\hat{s}_i$ does not
intersect the middle of the word $c_{i-1} u_{i-1} d_{i-1}$. But
$\hat{s}_i$ could not be situated in the first half of this word
because otherwise the first half of the word $q_{i+1}$ would be
strictly less than $u$ and consequently $q_n<uu^*$ which is a
contradiction. Thus $\hat{s}_i$ is a subword in the right half of
the word $q_{j+1}$ unlike that of $q_n=uu^*$ is intersected by the
unshrinkable word $u_j$. Thus there is $k\in \{j+1,\ldots, n-1\}$
such that $u_j$ and $\hat{s}_k$ do form  a  composition. This
contradicts to the conditions of the theorem.
\end{proof}

1. Let $S= \{w_j\}_{j\in \Re}$ be a symmetric set consisting of
unsrinkable words. Since compositions of any two words are always
zero  this set is closed under compositions. The other conditions
in the definition of strictly non-expanding set is obvious. Thus
$*$-algebra
\[
\mathbb{C}\big\langle x_1,\ldots, x_n, x_1^*,\ldots,x_n^* | w_j,\
j\in \Re \big\rangle
\]
 has a faithful $*$-representation by unbounded
operators.

2.  Let us consider in more detail the simplest example of
monomial $*$-algebras $\mathfrak{A}_{x^2}=\algebra{x,x^*|x^2=0,
x^{*2}=0}$. It was proved in~\cite{tapper} that $*$-algebra
$\algebra{x,x^*|x^p=0, x^{*p}=0}$ is $C^*$-representable for every
integer $p$.  We will show  that among the representations given
by theorem~\ref{main} there is a $*$-representation in bounded
operators. We believe that this is true for any monomial
$*$-algebra with unshrinkable relations but even for
$\algebra{x,x^*|x^3=0, x^{*3}=0}$ we could not make an explicit
calculations as we do for $\mathfrak{A}_{x^2}$.

It can be easily verified that $BW$ consists of words
$u_k=x(x^*x)^k$, $v_k=x^*(xx^*)^k$, $a_m=(xx^*)^m$, $b_m=(x^*x)^m$
where $k\ge 0, m\ge 1$. Obviously $BW_+$ consists  of only words
$a_m$ and  $b_m$  ($m\ge 1$). If $z$ ,$w\in BW$ then $zw^*\in W_+$
only in the following cases

\noindent 1.)  $z=u_k, w=u_t$; 2.)  $z=v_k, w=v_t$; 3.) $z=a_m,
w=a_n$; 4.)  $z=b_m, w=b_n$.
 Moreover,
$$ u_ku^*_t=a_{k+t+1},
v_kv^*_t=b_{k+t+1}, a_m a^*_n=a_{n+m},b_m b^*_n=b_{n+m}.$$
Consider the following ordering $$u_0<u_1<\ldots < a_1<a_2<\ldots
<v_0<v_1<\ldots < b_1<b_2<\ldots .$$ Denote
$\alpha(a_m)=\alpha_m$, $\alpha(b_m)=\beta_m$ then the Gram matrix
of the inner product defined in theorem~\ref{scal} is block
diagonal
$$\begin{pmatrix}
  A & 0 & 0 & 0 \\
  0 & A' & 0 & 0 \\
  0 & 0 & B & 0 \\
  0 & 0 & 0 & B'
\end{pmatrix}.$$ Where $A, A', B, B'$ are Gankel matrices  $A=(a_{ij})_{ij}$, $a_{ij}=\alpha_{i+j-1}$,
$A'=(a'_{ij})_{ij}$, $a'_{ij}=\alpha_{i+j}$, $B=(b_{ij})_{ij}$,
$b_{ij}=\beta_{i+j-1}$, $B'=(b'_{ij})_{ij}$,
$b'_{ij}=\beta_{i+j}$. Note that $Y'$ obtained from $Y$ by
cancelling out the first column (here $Y$ stands for $A$ or $B$).

Thus the question of positivity of the form $\sca{\cdot}{\cdot}$
is reduced to the question of simultaneous positivity of two
Gankel matrices $C$ and $C'$ where the second is obtained by
cancelling out the first column. We will show that such matrices
$A, A', B, B'$ could be chosen to be positive and such that $B=A$
and that the representation in theorem~\ref{main} is in bounded
operators.

Let $f:[0,1]\to [0,1]$ be continuous function $f(x)>0$ for all
$x\in [0,1]$. Let \[ \alpha_m = \int^1_0 t^{m+1} f(t) dt \] be the
moments of the measure with density $f(t)$.  It is well known that
then moment matrix $A=(a_{ij})_{i,j=1}^n$ (where $a_{ij}=
\alpha_{i+j-1}$) is positively defined. But then $A'$ is the
moment matrix of the measure with density $t f(t)$ and thus is
also positively-defined. We can put $B=A$.

To prove that the representation is in bounded operators we need
only to verify that the multiplication $L_x$ by  generator $x$ is a bounded
operator. Obviously, $x u_k=0$ and $xa_m =0$  for all $k\ge 0$ and
$ m\ge 1$. Compute $||xv_k||^2= \sca{a_{k+1}}{a_{k+1}}=
\alpha_{2(k+1)}$, $ ||v_k||^2=
\alpha(b_{2k+1})=\beta_{2k+1}=\alpha_{2k+1}$. Analogously,
$||xb_k||^2=\alpha_{2k+1}$ and $||b_k||^2 = \alpha_{2k}$. Thus
$L_x$ is bounded if there is a constant $c\ge 0 $ such that
\[
\alpha_{2(k+1)}\le c  \alpha_{2k-1}, \quad
 \alpha_{2k+1}\le
c\alpha_{2k}
\] for all $k\ge 1$. But
\[
\alpha_{2k}= \int^1_0 t^{2k+1} f(t) dt \le \int^1_0 t^{2k} f(t) dt
= \alpha_{2k-1}
\] and
\[
\alpha_{2k+1}= \int^1_0 t^{2k+2} f(t) dt \le \int^1_0 t^{2k+1}
f(t) dt = \alpha_{2k}.
\] Thus $||L_x||\le 1$. This proves that $\mathfrak{A}_{x^2}$ is $C^*$-representable.

3.  The $*$-algebra given by the generators and relations:
\[
\mathbb{C}\big<a_k,k=1,n | a_i^*a_j=\sum_{k\neq l}
T_{ij}^{kl}a_la_k^*;i\neq j  \big>,
\]
 with matrix of coefficient
satisfying $T_{ij}^{kl}=\bar{T}_{ji}^{lk}$ is strictly appropriate
by corollary~\ref{cor}.  since no two elements from defining
relations form a composition and the greatest word of any relation
begins with some $a_j$ and all other words begin with some
$a_k^*$.
 Hence it has a faithful unbounded representation. Note that with additional relations
$a_i^*a_i=1 +\sum_{k,l} T_{ii}^{kl}a_la_k^*$ we obtain so called
Wick's algebras.

4. If  $S \subset \mathbb{C}W(x_1,\ldots,x_n)$  is closed under
 compositions  then the $*$-algebra

\[
A=\algebra{x_1,\ldots,x_n,x^*_1,\ldots,x^*_n|\  S\cup S^*}
\]
(sometimes called $*$-double algebra)  is non-expanding (see
corollary~\ref{dub}).

If $S$ satisfies additionally the  property that  the greatest
word of any relation begins with generator different from the
beginnings of other longest words of this relation then $A$ is
strictly appropriate since $S\cup S^*$ is closed under
compositions. In particular, if $B$ is finite dimensional
associative algebra then its table of multiplication form a set of
relations $S$ with the greatest words of length $2$ and others of
length $1$. Thus $*$-algebra $A$ which $*$-double of $B$. In other
words,  $A$ is a free product $B_1*B_2$,  where $B_1\backsimeq B_2
\backsimeq B$ and involution is given on the generators
$b^*=\overline{\phi(b)}$ for any $b\in B_1 $ and
$c^*=\overline{\phi^{-1}(c)}$ for any $c\in B_2$ with $\phi:
B_1\to B_2$ being any isomorphism.

\begin{theorem}
Let $S\subset F$ be a symmetric subset of a free $*$-algebra in
generators $x_1,\ldots, x_n$ and $x^*_1,\ldots, x^*_n$ closed
under compositions such that for any $s\in S$ the following
properties holds:
\begin{enumerate}
\item $\hat{s}\in G=W(x_1,\ldots, x_n)$ or
$\hat{s}\in G^*=W(x^*_1,\ldots, x^*_n)$.
\item for any $u\prec s$ such that $|u|=|\hat{s}|$ words $u$ and
$\hat{s}$ both lie in the same semigroup $G$ or $G^*$.
\end{enumerate}

Then $S$ is non-expanding.
\end{theorem}
\begin{proof}
If some word $w\in W$ ($w=y_{k_1}\ldots y_{k_t}$ where $y_{k_r}\in
X\cup X^*$ are generators) contains subword $\hat{s}$ with $s\in
S$. Then  $w=p\hat{s}q$ for some words $p$ and $q$ in $W$. After
substitution $w\to p\bar{s}q=\sum_i \alpha_i p w_i q$ we see that
all words $w_i$ such that $|w_i|=|\hat{s}|$ lie in the same
semigroup either in $G$ or in $G^*$. Since decomposition ${\rm
R}_{S}(w)=\sum_j \beta_j u_j, u_j\in BW, u_j=z_{i_1}\ldots
z_{i_k}, z_{i_r}\in X\cup X^*$ can be obtained by several steps of
substitution considered above we see that for any $j$ such that
$|u_j|=|w|$ for all $r$ both generators $z_{i_r}$ and $y_{k_r}$
lie in the same set $X$ or $X^*$.

Let $u,v\in BW,\ u\not=v$. Let word $w$ be such that  $ ww^*\prec
{\rm R}_{S}(uv^*)$. If $|u|\not=|v|$ then $|w|<max(|u|,|v|)$ hence
$w<\sup{(u,v)}$. So assume that $|u|=|v|$ then $ww^*\le uv^*$
implies that $w\le u$. Let $w=u$. Write down $u=z_1\ldots z_k,
z_r\in X\cup X^*$. Without loss of generality assume $z_k\in X$.
Then $uu^*=z_1\ldots z_k z_k^*\ldots z_1^* $  by the first part of
the proof $u$ ends with generator from $X$ whereas $v^*$ begins
with generator from $X^*$. Thus $uv^*\in BW$ and ${\rm
R}_{S}(uv^*)=uv^*$. Clearly, $uv^*=uu^*$ implies $u=v$. Obtained
contradiction proves that if $ww^*\prec {\rm R}_{S}(uv^*)$ then
$w< \sup{(u,v)}$. Since  $S$ is appropriate for any  $d\in BW$
word $dd^*$ lies in $BW$ by lemma~2~\cite{pop1}.
\end{proof}

It could be shown using Zorn's lemma that for any algebra $A$ and
any its set of generators $S$ there is a Gr\"obner basis
corresponding to $S$ with any given inductive ordering of the
generators. Thus we have the following corollary.
\begin{corollary}\label{dub}
If $B$ is a finitely generated  associative algebra then its
$*$-double  $A=B*B$ is strictly non-expanding $*$-algebra. Hence
$A$ has a faithful $*$-representation in pre-Hilbert space.
\end{corollary}

Below we give some known examples of *-doubles which have finite
Gr\"obner bases.

5. Consider $*$-algebra:

\[Q_{4,\alpha}=\mathbb{C} \big\langle
q_1,\ldots,q_4,q^*_1,\ldots , q^{*}_4|\ q_j^2=q_j, \ \sum_j
q_j=\alpha, q_j^{*2}=q_j^*, \sum_j q_j^{*}= \overline{\alpha }
\big\rangle.
\] It is a *-double of the
algebra

\[
B_{n,\alpha}=\algebra{q_1,\ldots,q_4|\ q_j^2=q_j, \ \sum_j q_j=
\alpha }
\]
which has the following Gr\"obner basis:

$S=\{ q_1  q_1 - q_1, q_2 q_2 - q_2,
  -q_3 q_2 -2q_1 - 2q_2 - 2 q_3 + \alpha + 2 \alpha q_1+
      2 \alpha q_2 + 2  \alpha q_3 - \alpha^2 - q_1 q_2 - q_1 q_3 -
      q_2q_1 - q_2q_3 - q_3 q_1, q_3 q_3 -q_3,
  -q_3 q_1 q_2  -3\alpha + 5\alpha^2 - 2\alpha^3 +
      q_2(6 - 10\alpha + 4\alpha^2) +
      q_3(6 - 10\alpha + 4\alpha^2) +
      q_1(8 - 13\alpha + 5\alpha^2) + (3 -
            2\alpha)q_1  q_2 + (6 - 4\alpha) q_1 q_3 +
            (6 - 4\alpha)q_2 q_1 + (6 - 4\alpha) q_2 q_3 +
            (3 - 2\alpha)q_3 q_1 + q_1 q_2 q_1 +
      q_1 q_2q_3 + q_1q_3 q_1 + q_2q_1q_3 + q_2 q_3 q_1)
\} $.  More about this algebra can be found in ~\cite{sam,bart}.
Let us note that when $\alpha =0$ $*$-algebra
$Q_{4,0}=B_{4,0}*B_{4,0}$ has only zero representation in bounded
operators (see~\cite{bart}). Thus for this $*$-algebra unboundedness is essential.

6. That the generators in the previous example are idempotents is
not important,  we can consider other powers  as well:
\[T_{3,\alpha}=\mathbb{C} \big\langle q_1,q_2,q_3,q^*_1,q_2^* ,
q^{*}_3|\ q_j^3=q_j, \ \sum_j q_j=\alpha, q_j^{*3}=q_j^*, \sum_j
q_j^{*}= \overline{\alpha } \big\rangle.
\]
It is a *-double  of the algebra $\algebra{q_1,q_2,q_3|\
q_j^3=q_j, \ \sum_j q_j= \alpha }$. We can find its Gr\"obner
basis. We have the following set of relations
$
 \{q_1^3 - q_1,
 q_2^3 - q_2,
 q_3^3 - q_3,
 q_1 + q_2 + q_3 - \alpha
 \}
$. From these relations we see that this  algebra is generated by
$q_1$ and $q_2$. That is why we can consider the following set of
relations: $ \{q_1^3 - q_1,
 q_2^3 - q_2,
(\alpha - q_1 - q_2)^3-(\alpha - q_1 - q_2) \}$.  Consider the
following  order on generators $q_2>q_1$. All relations are
already  normalized. The greatest words in this relations are
$q_1^3$, $q_2^3$ and $q_1^2 q_2$. Thus we have no reductions. The
first and third relations give  two compositions. From one side
they intersect by the word $q_1$. And the result of this
composition is
  $(q_1^3 - q_1)q_1q_2-q_1^2((\alpha - q_1 - q_2)^3-(\alpha - q_1 - q_2))$.
  From other side  they  intersect by the word $q_1^2$. Result of this composition is
  $(q_1^3 - q_1)q_2-q_1((\alpha - q_1 - q_2)^3-(\alpha - q_1 - q_2))$.
  Another composition is formed by third and second relations.
  Their greatest words intersect by the word $q_2$. Result of this  composition is
  $((\alpha - q_1 - q_2)^3-(\alpha - q_1 - q_2))q_2^2- q_1^2(q_2^3 -
  q_2)$.
Hence we have three new relations. After performing reductions we
will have the following set of relations:

$
S=\{ q_1^3 - q_1,
  -q_2^2 q_1 + 3\alpha q_1^2 + 3\alpha q_2^2 + \alpha^3 +
      q_1(-1 - 3\alpha^2) + q_2(-1 - 3\alpha^2) +
      3\alpha q_1  q_2 - q_1 q_2^2 - q_1^2  q_2 + 3\alpha q_2 q_1 -
      q_2 q_1^2 - q_1 q_2 q_1 - q_2 q_1 q_2, q_2^3 - q_2,
  -q_2 q_1  q_2  q_1^2 + -\alpha^3 + 9\alpha^5 -
      q_1^2(-3\alpha - 37\alpha^3) -
      q_2^2(3\alpha - 27\alpha^3) -
      q_2(-1 + 6\alpha^2 + 27\alpha^4) -
      q_1(18\alpha^2 + 30\alpha^4) - (-12\alpha -
            45\alpha^3)q_1 q_2 -
      27\alpha^2 q_1q_2^2 -( 1 + 30\alpha^2) q_1^2  q_2 +
      9\alpha q_1^2 q_2^2 - (6\alpha - 18 \alpha^3) q_2q_1
      - (1 + 3\alpha^2) q_2 q_1^2 -
      (-2 + 15\alpha^2) q_1q_2 q_1 +
      3\alpha q_1q_2 q_1^2 + 3\alpha q_1^2 q_2 q_1 -
      q_1^2q_2 q_1^2 - (-1 + 9\alpha^2) q_2 q_1  q_2 +
      6\alpha q_1q_2 q_1q_2 - q_1^2 q_2 q_1 q_2 -
      3\alpha q_2 q_1q_2q_1 + q_1q_2q_1q_2q_1  \}
$

Some of these relations do form compositions but all of them
reduce to zero. Hence it is a Gr\"obner basis.

\section{APPENDIX: Noncommutative Gr\"obner bases.}
For the convenience of the reader we repeat the relevant material
from  noncommutative Gr\"obner bases theory (see~\cite{ufn,Bok})
with some easy reformulations.
 Let $W_n$ denote the  free $*$-semigroup  with   generators
 $x_1,\ldots, x_n$. For a word  $w= x_{i_1}^{\alpha_1} \ldots
 x_{i_k}^{\alpha_k}$ (where $i_1$, $i_2$, $\ldots$, $i_k\in \{1,\ldots, n\}$, and $\alpha_1$,
 $\ldots$, $\alpha_k \in \mathbb{N}\cup \{0\}$) the length of $w$,
 denoted by $|w|$, is defined to be $\alpha_1+\ldots +\alpha_k$.
  Let $F_n = \mathbb{C}\langle x_1, \ldots, x_n\rangle$ denote the free associative
algebra with   generators  $x_1,\ldots, x_n$. We will sometimes
omit subscript $n$.  Fix the  linear order on $W_n$ such that $x_1
> x_2 > \ldots > x_n$, the words of the same length ordered
lexicographically and the words of greater length are considered
greater. Any  $f \in F_n$  is a linear combination $\sum_{i=1}^k
\alpha_k w_i$ (we may assume that $\alpha_i\not=0$ for all $i\in
\{1,\ldots, k \}$) of words $w_1$, $w_2$, $\ldots$, $w_k$.  Let
$\hat{f}$ denote  the greatest of these words, say $w_j$. Then
denote $\hat{f}- (\alpha_j)^{-1}f$ by $\bar{f}$. The degree of
$f\in F_n$, denoted by $\deg(f)$, is defined to be $|\hat{f}|$.
Elements of the free algebra $F$ can be identified with functions
$f:W\to \mathbb{C}$ with finite support via the map $f\to
\sum_{w\in W}f(w)w$. For a word $z\in W$ and an element $f\in F$
we will write $z\prec f$ if and only if $f(z)\not= 0$.
\begin{definition}
We will say that two elements $f,g \in F_n$ form a composition
$w\in W$ if there are words $x,z \in W$ and nonempty word $y\in W$
 such that $\hat{f} = x y$, $ \hat{g} = y z$ and $w=x y z$,
 in other words  $w=\hat{f}w_1=w_2\hat{g}$ for some nonempty
 words $w_1$, $w_2\in W$ in which the marked words $\hat{f}$
 and $\hat{g}$ "intersect"  $|\hat{f}|>|w_2|$. Let us denote the result of the composition $\beta fw_1-\alpha w_2g$ by
 $(f,g)_w$,  where $\alpha$ and $\beta$ are the coefficients of the greatest words
  $\hat{f}$ and $\hat{g}$ in $f$ and $g$ respectively.
\end{definition}
It is obvious that $(f,g)_w < w$. Notice  that two elements $f$
and $g$ may form compositions in many ways and $f$ may form
composition with itself. The following definition is due to
Bokut~\cite{Bok}.

\begin{definition}
A subset $S \subseteq F_n$ is called closed under compositions  if
for any two elements  $f$, $g \in S$ the following properties
holds:
\begin{enumerate}
\item if $f\not=g$ then the word $\hat{f}$ is not a subword in
$\hat{g}$; \item  If $f$ and $g$ form a composition $w$ then there
are $a_j, b_j \in W_n,
 f_j \in S,  \alpha_j\in \mathbb{C}$ such that $(f,g)_w=\sum_{j=1}^m \alpha_ja_jf_jb_j$ and
 $a_jf_jb_j < w$, for $j=1,\ldots,m$.
\end{enumerate}
\end{definition}

\begin{definition}
A set $S\subseteq F$  is called a Gr\"obner basis of an  ideal
$\mathcal{I}\subseteq F$ if for any $f\in \mathcal{I}$ there is
$s\in S$ such that $\hat{s}$ is a subword in $\hat{f}$. A
Gr\"obner basis $S$ of $\mathcal{I}$ is called minimal if no
proper subset of $S$  is a Gr\"obenr basis of $\mathcal{I}$.
\end{definition}

If $S$ is closed under compositions then $S$ is a minimal
Gr\"obner basis for the ideal $\mathcal{I}$ generated by $S$
\cite{Bok}. Henceforth we will consider only minimal Gr\"obner
bases.  Thus we will say that $S$ is a Gr\"obner basis of an
associative algebra $A=F/\mathcal{I}$ if $S$ is closed under
composition and $S$ generate  $\mathcal{I}$ as an ideal of $F$.
Let $GB$ be a Gr\"obner basis for $A$, $\hat{GB}$ be the set of
greatest words $\hat{s}$ of elements of $s\in GB$  and $BW(GB)$ be
the subset of those words in $W_n$ that contain no word from
$\hat{GB}$ as a subword. It is a well known that $BW(GB)$ is a
linear basis for $A$. Further on we will write simply $BW$ since
we will always deal with a fixed Gr\"obner basis.

Let $S\subseteq F$ be closed under compositions, $\mathcal{I}$ the
ideal generated by $S$. Let us define the  operator ${\rm
R}_{S}:F\to F$ by the following rule. For $f\in F$  there are
uniquely defined coefficients $\{\alpha_i\}\subset \mathbb{C}$ and
words $\{w_i\}\subset BW$ such that
$f+\mathcal{I}=\sum_{i}\alpha_i(w_i+\mathcal{I})$ (basis
decomposition in the factor algebra). Put ${\rm
R}_{S}(f)=\sum_{i}\alpha_i w_i$. Then ${\rm R}_{S}$ is a
retraction on a subspace $K$ in $F$ spanned by $BW$. We can
consider the space $K$ with the new operation $f\diamond g= {\rm
R}_{S}(f g)$ for $f$, $g \in K$. Then $(K,+,\diamond)$ becomes an
algebra isomorphic to $F/\mathcal{I}$. The decomposition of
element $f=\sum_{j}\beta_j u_j$ where $u_j\in W$ with respect to
the basis $BW$ in the factor algebra $A/\mathcal{I}$ can be
obtained by repeated application of the following procedure: if
word $u_j$ contains subword $\hat{s}$ with $s\in S$, then  there
are words $p$ and $q$ such that $u_j=p\hat{s}q$.  Substitute
$\bar{s}$ instead of $\hat{s}$ (we will denote this substitution
by $p\hat{s}q\to p\bar{s}q$). Having done this for all $j$ we
obtain an element $\sum_{j}\beta_j^{'} u_j^{'}\in F$. Repeat this
procedure for all $u_j'$ and so on. After finite steps we obtain
desired decomposition $\sum_{i}\alpha_iw_i$ where all $w_j\in BW$.
From this follows that if $w\prec R_{S}(u)$ for some words $w$ and
$u$ then $w<u$.

\noindent {\bf Acknowledgements.}

 The authors wish to express their thanks to their
teacher Prof. Yurii Samoilenko for his kind attention and many
helpful comments.

\end{document}